\documentclass[12pt]{amsart}
\input xy
\xyoption{all}
\newcommand{\bP}{{\mathbb P}}
\newcommand{\bQ}{{\mathbb Q}}

\newcommand{\cD}{{\mathcal D}}

\newcommand{\cH}{{\mathcal H}}

\newcommand{\cO}{{\mathcal O}}

\newcommand{\ox}{\overline{x}}
\newcommand{\oz}{\overline{z}}

\newcommand{\tC}{\widetilde{C}}
\newcommand{\tD}{\widetilde{D}}

\newcommand{\tX}{\widetilde{X}}
\newcommand{\tY}{\widetilde{Y}}

\newcommand{\ra}{\rightarrow}
\newcommand{\lra}{\longrightarrow}

\newcommand{\inj}{\hookrightarrow}

\newcommand{\T}{\Theta}

\DeclareMathOperator{\tr}{tr}

\newtheorem{theoremi}{Theorem}
\newtheorem{theorem}{Theorem}[section]
\newtheorem{lemma}[theorem]{Lemma}
\newtheorem{proposition}[theorem]{Proposition}
\newtheorem{corollary}[theorem]{Corollary}
\newtheorem{remark}[theorem]{Remark}
\newtheorem{definition}[theorem]{Definition}
\newtheorem{hypothesis}[theorem]{Hypothesis}

\numberwithin{equation}{section}

\begin{document}

\title{Correspondences with split polynomial equations}

\author{E. Izadi}

\address{Department of Mathematics, University of Georgia, Athens, GA
30602-7403, USA}

\email{izadi@math.uga.edu}

\author{H. Lange}

\address{Mathematisches Institut, Bismarckstr. 1 1/2, 91054 Erlangen, Germany}

\email{lange@mi.uni-erlangen.de}

\author{V. Strehl}

\address{Institut f\"ur Informatik, Martensstr. 1, 91098 Erlangen, Germany}

\email{strehl@informatik.uni-erlangen.de} 

\thanks{This work was done while E. Izadi was visiting the
Mathematisches Institut in Erlangen at the invitation of H. Lange and
as part of the Faculty Exchange Program between the University of
Georgia and the University of Erlangen. E. Izadi is greatly indebted to
H. Lange for his generous hospitality and to the Universities of
Georgia and Erlangen for their support.}


\begin{abstract}

We introduce endomorphisms of special jacobians and show that they
satisfy polynomial equations with all integer roots which we compute. The
eigen-abelian varieties for these endomorphisms are generalizations of
Prym-Tyurin varieties and naturally contain special curves
representing cohomology classes which are not expected to be
represented by curves in generic abelian varieties.

\end{abstract}

\maketitle

\section*{Introduction}

\noindent

Let $A$ be an abelian variety of dimension $g$ and $\T$ a divisor on
$A$ representing a principal polarization. The minimal cohomology
class for curves in $A$ is
\[
\frac{ [\T ]^{ g-1 } }{ (g-1)! }.
\]
By a well-known result of Matsusaka \cite{matsusaka59} the minimal
class is represented by a curve $C$ in $A$ if and only if $(A ,\T)$ is
the polarized jacobian of $C$. Welters \cite{welters871} classified
the abelian varieties in which twice the minimal class is represented
by a curve. More generally, Prym-Tyurin varieties of index $m$ contain
curves representing $m$ times the minimal class.

A Prym-Tyurin variety $P$ of index $m$ is, by definition, the image of
$D-id$ in the jacobian $JC$ of a curve $C$ where $D$ is an
endomorphism of $JC$ satisfying the equation $(D-id)(D+(m-1)id) =
0$. The image of an Abel embedding of $C$ in $JC$ by the map $(D-id) :
JC\ra P$ is a curve representing $m$ times the minimal class in $P$
\cite{welters871}.

There are few explicit constructions of Prym-Tyurin varieties in the
literature.

In this paper we consider the more general situation where the
jacobian of a curve admits endomorphisms satisfying polynomials of
higher degree that can be decomposed into products of linear factors
with integer coefficients which we compute. So our endomorphisms have
integer eigen-values and, after isogeny, the jacobians of our curves
split into the product of the eigen-abelian varieties of the
endomorphism. The images of Abel embeddings of our curves will, after
isogeny (to obtain principally polarized abelian varieties), give
curves representing multiples of the minimal class in the
eigen-abelian varieties. In a future paper, we will compute the multiples 
of the minimal class that one obtains.
As in the case of Prym-Tyurin
varieties, these multiples will be
computable from the coefficients of the polynomial equations of the
endomorphisms.

The curves that we consider are immediate generalizations of
constructions of Recillas, Donagi and Beauville (see
e.g. \cite{recillas74}, \cite{donagi92}, \cite{beauville82}). Roughly
speaking, they are defined as follows (for details see Section
1). Suppose given a ramified covering $\rho_n: X \ra Y$ of degree $n$
of smooth projective curves and an
\'etale double cover $\tX \ra X$. Then a covering $\tC \ra
Y$ of degree $2^n$ can be defined as the curve parametrizing the
liftings of fibres of $\rho_n$ to $\tX$. Moreover, the involution on
$\tX$ induces an involution $\sigma$ on the curve $\tC$. Assuming the
ramification of $\rho_n$ is simple, we show that the curve $\tC$ is
smooth and that it has either one or two connected components. We
concentrate on the case where $\tC$ consists of two smooth connected
components $\tC_1$ and $\tC_2$. The computations in the case where
$\tC$ is irreducible yield polynomial equations similar to those
obtained for the case $n$ odd below. We shall not address this case in
this paper.

To be more precise, suppose first that $n = 2k + 1 \geq 3$. In this
case $\sigma$ induces an isomorphism $\tC_1 \ra \tC_2$ and we denote
$C = \tC_1$. Using the involution on $\tX$ we introduce a
correspondence $D$ on $C$.  Our first result is Theorem 2.4, which
says that $D$ satisfies an equation of degree $k$, integral over the
integers, whose coefficients are given by explicit recursion
relations. Denoting the induced endomorphism of the jacobian by the
same letter, clearly any integer zero of this equation yields an
eigen-abelian subvariety of $D$ on $JC$. Our main result for odd $n$
is that all zeros of this equation are integers. In fact, we have

\begin{theoremi} \label{mainthm1}
Suppose $n=2k+1$, $k\geq 1$. The correspondence $D$ satisfies the equation
\[
\prod_{ i=0 }^k (X + (-1)^{ i+k+1 } ( 2i+1 )) =0.
\]
which obviously does not have any multiple root.
\end{theoremi}

Suppose now $n = 2k \geq 2$. Then the involution $\sigma$ induces an
involution on each component $\tC^i$ for $i =1$ and 2, which we denote
by the same letter. Hence $J\tC_i$ decomposes up to isogeny into
the product of the Prym variety $P_i^{\sigma} := im (\sigma -id)$ of
$\sigma$ and its complement $B_i^{\sigma} : = im (\sigma + id )$. In
this case we introduce a correspondence $\tD_i$ on the curve $\tC_i$
which for $n \geq 6$ decomposes the abelian varieties $B_i^{\sigma}$
and $P_i^{\sigma}$ further. Again we compute the equation for the
correspondence $\tD_i$. This is a polynomial equation in $\tD_i$ and
$\sigma \tD_i$. Setting $\sigma = 1$, respectively $\sigma = -1$, we
obtain an equation for the endomorphism induced on $B_i^{\sigma}$,
respectively $P_i^{\sigma}$, the coefficients of which are given by
explicit recursion relations (see Theorems \ref{theorem3.6} and
\ref{theorem3.7}).  Again we prove that all zeros of these equations
are integers and thus lead to decompositions of the abelian varieties
$B_i^{\sigma}$ and $P_i^{\sigma}$ for $n \geq 6$ into eigen-abelian
subvarieties.
 
\begin{theoremi} \label{mainthm2}
{\em (1)} Suppose $n=4k$ with $k\geq 1$. For $i=1$ and 2 the
correspondence $\tD_i$ induces endomorphisms on $B_i^{\sigma }$ and
$P_i^{\sigma}$ satisfying the equations
\begin{itemize}
\item on $B_i^{\sigma }$: \hskip70pt $\prod_{j=0}^k (X -
8(k-j)^2 + 2k) = 0$,

\item on $P_i^{\sigma}$: \qquad \qquad  $\prod_{j=0}^{k-1} (X +
8(k-j)^2 - 10k +8 j +2) = 0$.
\end{itemize}
{\em (2)}  Suppose $n=4k-2$ with $k\geq 2$. For $i=1$ and 2 the
correspondence $\tD_i$ induces endomorphisms
on $B_i^{\sigma }$ and $P_i^{\sigma}$ satisfying the equations   
\begin{itemize}
\item on $B_i^{\sigma }$: \hskip70pt $\prod_{j=0}^{k-1} (X -
8(k-j)^2 + 10k - 8j - 3) = 0$,
\item on $P_i^{\sigma}$: \hskip67pt $\prod_{j=0}^{k-1} (X + 8(k-j)^2
- 18 k + 16 j + 9) = 0$.
\end{itemize}
\end{theoremi}

It is easy to see that the polynomials involved do not have multiple roots.
The main idea of the proofs of Theorems \ref{mainthm1} and
\ref{mainthm2} is to identify the fibres of the coverings $f: C \ra Y$
and $f_i: \tC_i \ra Y$ with sub-vector spaces of the space of bit
vectors of length $n$. This gives an additional structure on the
fibres, namely that of a Hamming scheme,
as known from algebraic combinatorics and coding theory 
(not to be confused with a scheme in the algebro-geometric sense). 
Using this we
associate to $D$ and $\tD_i$ endomorphisms of vector spaces for which
we can explicitly determine the eigenvalues and eigenvectors. \\

The contents of the paper are as follows: In Section 1 we recall the
$n$-gonal construction. In Section 2 we introduce the correspondence
$D$ and compute its equation in the odd-degree case. Section 3
contains the analogous computations for even $n$.  In section 4 we
provide the combinatorial tools needed for the proofs of Theorems 1
and 2, which are given in Section 5.  In Section 6 we give a system of
equations for the dimensions of the eigen-abelian varieties involved.
We use these equations to compute these dimensions explicitly for $n
\leq 10$.  Finally, section 7 contains a combinatorial remark related
to the situation of Theorem 1 which is worth noting.

\section{The n-gonal construction}

\subsection{The set up}

Let $Y$ be a smooth curve of genus $g_Y$, $X$ a cover of degree $n$ of 
$Y$ of genus $g_X$ and $\tX$ an \'etale double cover of $X$ which is
NOT obtained by base change from a double cover of $Y$:
\[
\tX\stackrel{\kappa}{\lra} X\stackrel{\rho_n}{\lra }Y.
\]
Then $Y$ embeds into the symmetric power $X^{ (n) }$ via the map
sending a point $y$ of $Y$ to the divisor obtained as the sum of its
preimages in $X$. Let $\tC\subset\tX^n$ be the curve defined by the
fiber product diagram

\begin{eqnarray} \label{eq1}
\tC & \lra & \tX^{ (n) } \nonumber \\
\downarrow & & \downarrow \kappa^{ (n) } \\
Y & \lra & X^{ (n) }. \nonumber 
\end{eqnarray}

In other words, the curve $\tC$ parametrizes the liftings of points of
$Y$ to $\tX$.


\begin{lemma}
If $\rho_n$ is at most simply ramified, then the curve $\tC$ is smooth.
\end{lemma}

\begin{proof}


Since $\tC$ was defined by the fiber product diagram \eqref{eq1},
the tangent space to $\tC$ is the pull-back of the tangent space of
$Y$. Away from the branch points of $\rho_n$, the map $\tC\ra Y$ is
\'etale and hence $\tC$ is smooth. The ramification points of $\tC$
over $Y$ can be described as follows. Let $y\in Y$ be a branch point
of $\rho_n$. Let $\ox$ be the ramification point of $\rho_n$ above $y$
and let $\ox_1 ,\ldots ,\ox_{n-2}$ be the remaining (distinct) points
of $X$ above $y$. Then a point of $\tC$ above $y$ is a ramification
point if and only if it is of the form $x+x' +x_1+\ldots x_{n-2}\in
\tX^{ (n)}$ where $x$ and $x'$ are the two points of $\tX$ above $\ox$
and $x_i$ is a point of $\tX$ above $\ox_i$ for $i=1 ,\ldots
,n-2$. The tangent space to $\tX^{(n)}$ at $x+x' +x_1+\ldots
x_{n-2}\in \tX^{ (n)}$ can be canonically identified with
\[
\cO_x (x) \oplus \cO_{x'} (x') \oplus_{i=1}^{n-2}\cO_{x_i }(x_i)
\]
and the tangent space to $X^{(n)}$ at $2\ox+\ox_1+ \ldots + \ox_{n-2}$
can be canonically identified with
\[
\cO_{2\ox} (2\ox) \oplus_{i=1}^{n-2}\cO_{\ox_i }(\ox_i).
\]
The differential of $\kappa^{(n)}$ sends $\cO_{x_i} (x_i)$
isomorphically to $\cO_{\ox_i }(\ox_i)$ and sends $\cO_x (x)$ and
$\cO_{x' }(x')$ both isomorphically to the subspace $\cO_{\ox}(\ox)$
of $\cO_{2\ox} (2\ox)$. Its kernel is therefore one-dimensional and it
follows that $\tC$ is smooth at $x+x' +x_1+\ldots x_{n-2}$ if and only
if the image of the tangent space $\cO_y(y)$ of $Y$ at $y$ is not
contained in the subspace
\[
\cO_{\ox} (\ox) \oplus_{i=1}^{n-2}\cO_{\ox_i }(\ox_i).
\]
Equivalently, if and only if the composite map
\begin{equation}\label{eqcomptang}
\cO_y(y) \lra \cO_{2\ox} (2\ox) \oplus_{i=1}^{n-2}\cO_{\ox_i }(\ox_i)
\lra \cO_{\ox} (2\ox)
\end{equation}
where the second map is the quotient by the image of the tangent space
of $\tX^{ (n)}$ is not zero.

Now choose a general map $Y \ra \bP^1$ of
degree $m$ with simple ramification disjoint from the branch locus of
$\rho_n$. Let $p\in\bP^1$ be the image of $y$ by this map. Define
$\widehat{C}$ by the pull-back diagram
\[
\begin{array}{ccc}
\widehat{C} & \lra & \tX^{ (mn) } \\
\downarrow & & \downarrow \kappa^{ (mn) } \\
\bP^1 & \lra & X^{ (mn) }.
\end{array}
\]
By \cite[8.13, a) p. 107]{welters81}, the curve $\widehat{C}$ is
singular exactly above the ramification of the map $Y \ra
\bP^1$. In particular, it is smooth above $p$. Applying our analysis
above to this case, this means that the composite map
\[
\cO_p (p)\lra \cO_{2\ox }(2\ox) \oplus_{ i=1}^{mn -2}\cO_{\ox_i
}(\ox_i) \lra \cO_{\ox }(2\ox)
\]
is an isomorphism. Here $\ox_i$, $i= n-2 ,\ldots ,mn-2$ are the other
points of of $X$ above p. It is now easy to see that after identifying
$\cO_p (p)$ with $\cO_y (y)$ via the differential of $Y\ra\bP^1$, this
map is equal to the map \eqref{eqcomptang} which is therefore also an
isomorphism.  This shows that $\tC$ is smooth.
\end{proof}

Now we investigate the number of connected components of the curve
$\tC$. We first have

\begin{lemma}

If $\rho_n$ is unramified, then the curve $\tC$ is a union of $2n$
disjoint copies of $Y$.

\end{lemma}

\begin{proof}
This follows
immediately from the fact that, locally, a small loop in $Y$ will lift
to $n$ disjoint loops in $X$ and $2n$ disjoint loops in $\tX$.
\end{proof}

So the case where $\rho_n$ is unramified is uninteresting from the
point of view of construction of abelian subvarieties of
jacobians. {\em From now on we will assume that $\rho_n$ is ramfied
with simple ramification.}

Recall that the Norm map $Nm : Pic^n \tX\ra Pic^n X$ is defined as
$\cO_{\tX } (D) \mapsto \cO_X (\kappa_* D)$ and that its kernel has
two connected components that are translates of the Prym variety $P$
of the double cover $\kappa :\tX\ra X$. Therefore the fibers of the
induced map
\[
Nm |_Y : Nm^{ -1 } (Y)\lra Y\subset X^{ (n) }\lra Pic^n (X)
\]
are disjoint unions of two translates of $P$. Let $\tY\ra Y$ be the
\'etale double cover parametrizing the components of the fibers of $Nm
|_Y$. Then, by the definition of $\tC$, the composite map
\[
\tC\inj \tX^{ (n) } \lra Pic^n \tX
\]
induces a map $\tC\ra \tY$ whose composition with $\tY\ra Y$ is the
natural map $\tC\ra Y$ from \eqref{eq1}. We have

\begin{lemma}

The curves $\tC$ and $\tY$ have the same number of connected
components.

\end{lemma}

\begin{proof}

As in the proof of \cite[Proposition 8.8 p. 100]{welters81} (also see
\cite[page 109]{I3}), it can be seen that any two points in a fiber of
$\tC\ra \tY$ can be joined by a path in $\tC$.
\end{proof}

From now on we make the following assumption.
\begin{hypothesis}\label{hyp}

The map $\rho_n$ is simply ramified and the double cover $\tY\ra Y$ is trivial.

\end{hypothesis}

By the above lemmas, this is equivalent to the fact that $\tC$ is smooth with two
connected components. Note that when $Y\cong\bP^1$, the double cover $\tY\ra Y$ is always trivial.

One situation (see \cite[Section 2.2]{donagi92}) in which $\tY\ra Y$
is trivial is when $\tX\ra X\ra Y$ is simple of type $D_n$, i.e., has
the following properties.

\begin{definition}
We say that the covering
$\tX\stackrel{\kappa}{\lra} X\stackrel{\rho_n}{\lra }Y$ is a {\it
simple covering of type} $D_n$ if
\renewcommand{\theenumi}{\roman{enumi}}
\begin{enumerate}
\item 
$\rho_n: X \ra Y$ is simply ramified of degree $n$ with branch divisor
$\cD \neq \emptyset$ and $\kappa: \tX \ra X$ an \'etale double
covering;
\item 
$\rho_n: X \ra Y$ is a primitive covering;
\item
the monodromy map of the covering $\rho_n \circ \kappa: \tX \ra Y$ can
be decomposed as
$$
\pi_1(Y \setminus \cD, y_0) \ra W(D_n) \hookrightarrow S_{2n}.
$$
\end{enumerate}
\end{definition}

Here $y_0 \in Y \setminus \cD$ is a base point, $W(D_n)$ denotes the
Weyl group of type $D_n$ and $W(D_n) \hookrightarrow S_{2n}$ the
standard embedding. Recall that a covering is called primitive if it
is not the composition of two coverings of degree $\geq 2$. The simply
ramified covering $\rho_n$ is primitive if and only if the canonical
map $\pi_1(X,*) \ra \pi_1(Y,*)$ is surjective. According to
\cite[Corollary 2.4]{donagi92}, any covering
$\tX\stackrel{\kappa}{\lra} X\stackrel{\rho_n}{\lra }\bP^1$ satisfying
(i) and (ii) is a simple covering of type $W(D_n)$.\\

In general the curve $\tC$ can be irreducible. For
examples see \cite{kl} and use, in particular, Remark 2.10.

The involution $\sigma$ exchanging complementary liftings
of the same point of $Y$ acts on $\tC$ and we let $C$ be the quotient
of $\tC$ by this involution. This means the following. Let $\oz
:=\ox_1 +\ldots + \ox_n$ be the sum of the points in a fiber of
$\rho_n$, and, for each $i$, let $x_i$ and $x_i'$ be the two preimages
of $\ox_i$ in $\tX$. Then
\[ 
z := x_1 +\ldots + x_n
\] 
is a point of $\tC$
and 
\[
\sigma (z) = x_1' +\ldots + x_n'.
\] 
The degrees of the maps
$\tC\ra Y$ and $C\ra Y$ are $2^n$ and $2^{ n-1 }$
respectively. Since the ramification of $\rho_n$ is simple, it
is easily seen that $\sigma$ is fixed-point-free if $n\geq 3$. Also,
we can see that for each ramification point $\ox_1 =\ox_2$ of $\rho_n$
there are $2^{n-2}$ ramification points in a fiber of $\tC\ra Y$
obtained as $x_1 + x_1' + D_{ n-2 }$ where $D_{ n-2 }$ is one of the
$2^{ n-2 }$ divisors on $\tX$ lifting $\ox_3 +\ldots +\ox_n$.

Let $\tC_1$ and $\tC_2$ be the two connected components of $\tC$. Then
half of the divisors $x_1 + x_1' + D_{ n-2 }$ lie in $\tC_1$ and the
other half lie in $\tC_2$.

Writing the degree of the ramification divisor of $\rho_n$
as
\[
deg (R_{ X/Y } ) = 2g_X -2 - n( 2g_Y -2 ),
\]
this shows that the genus of $\tC_1$ and $\tC_2$ is
\[
g_{\tC_i } = 2^{ n-3 } \left( g_X -1 - (n-4 )(g_Y -1 ) \right) +1.
\]

If $n$ is odd, the involution $\sigma$ exchanges the two components of 
$\tC$, hence induces isomorphisms
\[
\tC_1\cong\tC_2\cong C.
\]
So we have the following diagram
\begin{eqnarray}\label{diag1.3}
\xymatrix {
\tC_1\cup\tC_2 \ar[d]_{2:1} & \tX \ar[d]^{2:1} \\
C =\tC_1 =\tC_2 \ar[dr]_{2^{n-1}:1} & X \ar[d]^{n:1} \\
 & Y.}
\end{eqnarray}

If $n$ is even, the involution $\sigma$ acts on each component of $\tC$ 
hence $C$ also has two connected components, say $C_1$ and $C_2$. For 
$n\geq 4$, since $\sigma$ is fixed-point-free, we compute the genus of 
$C_1$ and $C_2$ to be
\[
g_{ C_i } = 2^{ n-4 } \left( g_X -1 - (n-4 )(g_Y -1 ) \right) +1.
\]
In this case we obtain the diagram
\[
\xymatrix{
\tC_1 \ar[d]_{2:1} & \tX \ar[d]^{2:1} & \tC_2 \ar[d]^{2:1}\\
C_1 \ar[dr]_{2^{n-2}:1} & X \ar[d]^{n:1} & C_2 \ar[dl]^{2^{n-2}:1}\\
 & Y.}
\]

If $n=2$, the degree of each component $C_i$ over $Y$ is $1$ so 
\[
C_1\cong C_2\cong Y.
\]

\subsection{Notation} For each $k\in\{ 0,\ldots
,n\}$, we denote by
\[
[k + (n-k )'](z)
\]
the sum of all the points where
$k$ of the $x_i$ are added to $(n-k)$ of the $x_i'$, the indices $i$
being all distinct. For instance
\[
[1+ (n-1)'](z) =\sum_{1\leq i \leq n} x_1' +\ldots +x_{ i-1 }' + x_i +
x_{ i+1 }' +\ldots + x_n'.
\]
and
\[
[2+ (n-2)'](z) =\sum_{1\leq i < j\leq n} x_1' +\ldots +x_{ i-1 }' + x_i +
x_{ i+1 }' +\ldots + x_{ j-1 }' + x_j + x_{ j+1 }' +\ldots + x_n'.
\]

\section{The correspondence for $n$ odd}\label{sectodd}

\subsection{Definition of $D$} For $i=1$ or $2$, we define a correspondence
$D_i$ on $\tC_i$ as the reduced curve
\[
D_i :=\{ (x_1 +\ldots + x_n ,x_1 + x_2' +\ldots + x_n'
)\}\subset\tC_i\times\tC_i
\]
and we define
\[
D\subset C\times C
\]
as the image of $D_i$ in $C\times C$. Note that the image of $D_1$ in
$C\times C$ is equal to the image of $D_2$. The correspondence $D$
defines an endomorphism of the jacobian $JC$ whose ``eigenspaces'' are
proper abelian subvarieties of $JC$. We call these the eigen-abelian
varieties of $D$. The aim of this section is to determine the
polynomial equation satisfied by this endomorphism. To study this
correspondence, we work on the curve $C$ which we consider as $C_1$.

For any $z = x_1 +\ldots + x_n\in C$ we define as usual
\[
D(z)
= {p_2}_* ((p_1^* z)\cdot D )
\]
as divisors on $C$, where $p_1$ and $p_2$ are
the first and second projections. The points of $C$ in the support
of $D (z)$ are sums of $x_i$ or $x_i'$. 
It is immediate that
\begin{equation}\label{eqnD1odd}
D (z) = [1+ (n-1 )'](z)
\end{equation}
and
\begin{equation}\label{eqnD2odd}
D^2 (z) = nz + 2[2' +(n-2 )](z)
\end{equation}
where $D^i$ is the composition of $D$ with itself $i$ times.

\subsection{The general equation for $n$ odd}
Applying $D$ to successive equations, we can
find polynomial equations for $D$ for any $n$ odd. First we have

\begin{proposition}\label{propDk}
{\em (1)} For any even integer $k,\; 0 \leq k\leq\frac{n-2}{2}$, there
are integers $a^\ell_j$ satisfying an equation
\begin{equation}\label{eqnkeven}
D^k (z) = a_0^k z + a_2^k D^2 (z) +\ldots + a_{ k-2 }^kD^{ k-2 }
(z) + k! [k' + (n-k )] (z)
\end{equation}
{\em (2)} For any odd integer $1 \leq k\leq\frac{n-2}{2}$, there are
integers $a^\ell_j$ satisfying an equation
\begin{equation}\label{eqnkodd}
D^k (z) = a_1^k D (z) + a_3^k D^3 (z) +\ldots + a_{ k-2
}^kD^{ k-2 } (z) + k! [k + (n-k )'] (z).
\end{equation}
\end{proposition}

Note that the integers $a^\ell_j$ are defined only for $\ell \equiv j
\;\mbox{mod} \;2$, $0 \leq j < \ell \leq \frac{n-2}{2}$.
\begin{proof} According to equations (\ref{eqnD1odd}) and
(\ref{eqnD2odd}) the proposition is valid for $k = 0, 1$ and $2$.
Applying $D$ to (\ref{eqnkeven}), we obtain
\[
\begin{split}
D^{ k+1 } (z) = a_0^k D (z) + a_2^k D^3(z) +\ldots +a_{ k-2 }^k
D^{ k-1 } (z) + k ! (n-k +1 ) [(k-1) + (n-k+1 )'] (z) \\
+ (k+1 )! [(k+1)+ (n-k-1)'] (z).
\end{split}
\]
Using (\ref{eqnkodd}) to substitute for $[(k-1)' + (n-k+1 )] (z)$, this
becomes
\begin{equation}\label{eqnk+1odd}
\begin{split}
D^{ k+1 } (z) = (a_0^k - k(n-k+1) a_1^{ k-1}) D (z) + (a_2^k -
k(n-k+1) a_3^{ k-1})D^3(z) + \\ \vdots \\+ (a_{ k-4 }^k -k (n-k+1) 
a_{ k-3 }^{ k-1 }) D^{ k-3 } (z) \\
+ (a_{ k-2 }^k + k( n-k+1 ))D^{ k-1 } (z) + (k+1 )! [(k+1) + 
(n-k-1)'] (z).
\end{split}
\end{equation}

Similarly, applying $D$ to (\ref{eqnkodd}) and using
(\ref{eqnkeven}) to substitute for $[(k-1)' + (n-k+1)] (z)$, we obtain
\begin{equation}\label{eqnk+1even}
\begin{split}
D^{ k+1 } (z) = -k (n-k+1 ) a_0^{ k-1 } z + (a_1^k -k( n-k+1 )
a_2^{k-1}) D^2 (z) + \\ \vdots \\+ (a_{k-4 }^k -k(n-k+1 ) a_{ k-3 }^{ 
k-1
})D^{ k-3 } (z) \\
+ (a_{ k-2 }^k + k(n-k+1 ))D^{ k-1 } (z) + (k+1 )!
[(k+1)' + (n-k-1) ] (z).
\end{split}
\end{equation}
By induction this completes the proof.
\end{proof}

The proof of the proposition gives the following recursion relations
for the integers $a^\ell_j$.

\begin{corollary}
Setting $a^\ell_{\ell} = -1$ for $0 \leq \ell \leq \frac{n-4}{2}$ and
$a^{-1}_1 = a^\ell_{-1} = 0$ for odd $\ell$, we have
for all $i \equiv k+1 \;\mbox{mod}\; 2$ and $0 \leq i \leq k-1, $
\begin{equation} \label{recurs} 
a^{k+1}_i = a^k_{i-1} - k(n-k+1)a^{k-1}_i.
\end{equation}
\end{corollary}

Using this we obtain

\begin{proposition} \label{prop2.3}

With the above notation we have
\[
a_{ k- 2i }^k = (-1)^{ i+1 }\sum_{ j_1=j_0 + 2 }^{ k-2i+1 } j_1 (n-j_1+1)
\cdot\sum_{ j_2=j_1+2 }^{ k-2i+3 } j_2 (n-j_2+1) \cdot\ldots
\cdot\sum_{ j_{ i+1 }= j_i+2 }^{ k-1 } j_{ i+1 } (n- j_{ i+1 } +1)
\]
for $0\leq i\leq\frac{k+1}{2}$ and $k \leq \frac{n-2}{2}$, where we
set $j_0 = -1$.

\end{proposition}

\begin{proof}

We prove the formula by induction. The formula holds trivially for
$k=0$ and $1$. Assume now that it holds for all $\ell\leq k-1$ and all
$i$, $0\leq i\leq\frac{\ell+1}{2}$. We need to prove it for $\ell=k$ and all 
$i$, $0\leq i\leq\frac{k+1}{2}$. From \eqref{recurs} we deduce
\begin{eqnarray*}
a_{ k-2i }^k = a_{ k-2i-1 }^{ k-1 } - (k-1)(n-k+2) a_{ k-2i }^{ k-2 }
\hspace{12cm}\\
= (-1)^{ i+1 } \sum_{ j_1=1 }^{ k-2i } j_1 (n-j_1+1)
\cdot\sum_{ j_2=j_1+2 }^{ k-2i+2 } j_2 (n-j_1+1) \cdot\ldots
\cdot\sum_{ j_{ i+1 }= j_i+2 }^{ k-2 } j_{ i+1 } (n- j_{ i+1 } +1)
\hspace{3.5cm}\\
+ (-1)^{ i+1 } (k-1) (n-k+2) \sum_{ j_1=1 }^{ k-2i+1 } j_1 (n-j_1+1)
\cdot\sum_{ j_2=j_1+2 }^{ k-2i+3 } j_2 (n-j_1+1) \cdot\ldots
\cdot\sum_{ j_{ i }= j_{ i-1 }+2 }^{ k-1 } j_{ i } (n- j_{ i } +1)\hspace{2cm}
\end{eqnarray*}

Note that in the expression in the proposition, if we remove the
last term, i.e., $(k-1)(n-k+2)$, all the upper bounds of the former
sums go down by $1$. This gives us the second line above. The rest
will then be the third line above, which proves the proposition.
\end{proof}

\subsection{The final equations for $n$ odd} With these coefficients
$a^\ell_j$ the following theorem gives the equation for
the correspondence $D$.

\begin{theorem} \label{thm2.4}
Suppose $n=2k+1$. \\
{\em (1)} For $k$ even $D$ satisfies the equation
\begin{equation} \label{eqn2.8}
X^k +(k+1)\sum_{i=0}^{k-2}a^k_i X^i_1 - \sum_{i=1}^{k-1}a^{k+1}_i X^i_1 = 0.
\end{equation}
{\em (2)} For $k$ odd $D$ satisfies the equation
\begin{equation} \label{eqn2.9}
X^k +(k+1)\sum_{i=0}^{k-1}a^k_i X^i_1 - \sum_{i=1}^{k}a^{k+1}_i X^i_1 = 0.
\end{equation}
\end{theorem}

\begin{proof}
(1) If $n = 2k+1$ with $k$ even, then
 $[(k+1) + (n-k-1)'] (z) = [k' + (n-k)] (z)$ and we can use 
(\ref{eqnkeven}) to substitute in equation (\ref{eqnk+1odd}) which then 
becomes
\[
\begin{split}
D^{ k+1 } (z) = -(k+1) a_0^k z + (a_0^k - k (k+2) a_1^{ k-1})D (z) 
+ \\
\vdots \\
-( k+1) a_{k-4 }^kD^{ k-4 } (z) + (a_{k-4 }^k - k (k+2 ) a_{ k-3 
}^{k-1} )D^{ k-3 } (z) \\
- (k+1 ) a_{ k-2 }^kD^{ k-2 } (z) + ( a_{ k-2 }^k + k (k+2 ))D^{ 
k-1 } (z) + (k+1)D^k (z).
\end{split}
\]
From the recursion relations (\ref{recurs}) we see that this equation is
just \eqref{eqn2.8}.

(2) If $n=2k+1$ with $k$ odd, then $[ (k+1)' + (n-k-1)] (z) = [k +
(n-k)'] (z)$ and we can use
(\ref{eqnkodd}) to substitute in equation (\ref{eqnk+1even}) which then 
becomes \eqref{eqn2.9}, again using the recursion relations (\ref{recurs}).
\end{proof}

\section{The correspondence $\tD_i$ for $n\geq 4$ even}

\subsection{Definition of $\tD_i$} For $i=1$ or $2$, we define a correspondence
$\tD_i$ on $\tC_i$ as the reduced curve
\[
\tD_i :=\{ (x_1 +\ldots + x_n ,x_1 + x_2 + x_3'+\ldots + x_n'
)\}\subset\tC_i\times\tC_i.
\]
For $n\geq 6$ the map from $\tD_i$ onto its image in $C_i\times C_i$
is of degree $2$ and we define
\[
D_i\subset C_i\times C_i
\]
as the reduced image of $\tD_i$ in $C_i\times C_i$. For $n=4$, the map 
from $\tD_i$ onto its image in $C_i\times C_i$ is of degree $4$ and we 
define
\[
D_i\subset C_i\times C_i
\]
to be twice the reduced image of $\tD_i$ in $C_i\times C_i$.

The correspondences 
$D_i$ and $\tD_i$
define endomorphisms of the jacobians $JC_i$ and $J\tC_i$
whose eigen-abelian varieties are proper abelian subvarieties of $JC_i$ and 
$J\tC_i$. The aim of this section is to determine the polynomial
equations satisfied by
these endomorphisms.

As before, for any $z = x_1 +\ldots + x_n\in\tC_i$ write
\[
\tD_i(z)
= {p_2}_* ((p_1^* z)\cdot\tD_i )
\]
as divisors on $\tC_i$, where $p_1$ and $p_2$ are
the first and second projections. With the notation of Section 1.2, we have
\begin{equation}\label{eqnD1even}
\tD_i (z) = [2+ (n-2 )'](z)
\end{equation}
and
\[
\tD^2_i (z) = { n\choose 2 }z + 2 (n-2) [2' +(n-2 )](z) + 6 [4' +
(n-4) ](z),
\]
which can be rewritten as
\begin{equation}\label{eqnD2even}
\tD^2_i (z) = { n\choose 2 }z + 2 (n-2)\sigma\tD_i (z) + 6 [4' +
(n-4) ](z).
\end{equation}

\begin{remark}
If $n=2$, then the correspondences $\tD_i$ are just the diagonals
of $\tC_1\times \tC_1$ and $\tC_2\times\tC_2$. So $D_1$ and $D_2$ are
the diagonals of $C_1$ and $C_2$.
\end{remark}

\subsection{Splitting of the jacobians}\label{ss3.2} The involution
$\sigma$ splits the jacobians of $\tC_i$
into their $+1$ and
$-1$ eigen-abelian varieties, i.e., the respective images of $\sigma +1$ and
$\sigma -1$. We denote
\[
P_i^{\sigma } := Im (\sigma -1 )\subset J\tC_i\qquad B_i^{\sigma } := 
Im
(\sigma +1 )\subset J\tC_i.
\]
Note that $B_i^{\sigma }$ is the image of $JC_i$ by the pull-back map 
of $\tC_i\ra C_i$.

It is immediate from the definitions that the
endomorphisms $\sigma$ and $\tD_i$ commute on $J\tC_i$. Hence $\tD_i$
induces endomorphisms on $P_i^{\sigma }$ and $B_i^{\sigma }$ which we
denote again by $\tD_i$.

As the double cover $\tC_i\ra C_i$ is \'etale, the map $J C_i\ra
B_i^{\sigma }$ which is obtained from pull-back of line bundles from
$C_i$ to $\tC_i$ has degree $2$. The endomorphism of $J C_i$ obtained
from $D_i$ and that of $B_i^{\sigma }$ obtained from $\tD_i$ fit into
the commutative diagram
\[
\xymatrix{
J C_i \ar[d] \ar[r]^{ D_i } & J C_i \ar[d] \\
B_i^{\sigma } \ar[r]^{ \tD_i } & B_i^{\sigma }.}
\]

\subsection{The general equation for $n$ even}
We proceed as in the case $n$ odd to find the general equation for
$\tD_i$, for $i=1$ or $2$. In order to formulate it, we define 
\[
\{ k\} :=\prod_{ i=1 }^k {2i\choose 2}.
\]
\begin{proposition}\label{propDkeven}
For $i = 1$ and 2 and any integer $k$, $2\leq k\leq\frac{n-2}{4}$,
there are integers $b_j^k$, $0 \leq j \leq k$ satisfying an equation
\begin{equation}\label{eqn3.3}
\tD_i^k (z) = \sum_{j=0}^{k-1}b_j^k\sigma^{k+j}\tD_i^j (z) 
+ \{ k\}\sigma^k [(2k)' + (n-2k )] (z).
\end{equation}
\end{proposition}
Note that $\sigma^\ell= id$ for $\ell$ even and $\sigma^\ell =\sigma$
for $\ell$ odd.
\begin{proof}
Suppose first $k=2$. Then
\begin{eqnarray*}
\tD^2_i (z) = &{ n\choose 2 }z + 2 (n-2) [2' +(n-2 )](z) + 6 [4' +(n-4) ](z)\\
  = & { n\choose 2 }z + 2 (n-2)\sigma\tD_i (z) + 6 [4' +(n-4) ](z).
\end{eqnarray*}
which is of the form (\ref{eqn3.3}). For $2\leq k\leq\frac{n-4}{4}$ we 
apply $\tD_i$ to (\ref{eqn3.3}) to obtain
\begin{equation}\label{eqnknevenrecur}
\begin{split}
\tD_i^{ k+1 } (z) = b_0^k \sigma^k\tD_i (z) + b_1^k\sigma^{ k+1 }
\tD_i^2(z) +\ldots +b_{ k-1 }^k\sigma^{ 2k-1 }
\tD_i^{ k } (z) \\ +
\{ k\}{ n-2k+2\choose 2}\sigma^{ k+1 } [(2k-2)' +
(n-2k+2)] (z) \\ + \{ k\} 2k (n-2k )\sigma^{ k+1 } [(2k)' + (n-2k)] (z)
\\ +\{ k\} {2k+2\choose 2}\sigma^{ k+1 } [(2k+2)' + (n-2k-2) ] (z).
\end{split}
\end{equation}
First assume $k\geq 3$. Then, using (\ref{eqn3.3}) to substitute
for $[(2k-2)' + (n-2k+2 )] (z)$ and $[ (2k)' + (n-2k)] (z)$, we can
write
\[
\begin{split}
\{ k\}{ n-2k+2\choose 2}\sigma^{ k+1 } [(2k-2)' +
(n-2k+2)] (z) \hspace{2cm}\\ = { 2k\choose 2 }{ n-2k +2\choose
2}\left(\tD_i^{ k-1 }
(z) -\sum_{ j=0 }^{ k-2 } b_j^{ k-1 }\sigma^{ k-1+j }\tD_i^j
(z)\right)
\end{split}
\]
and
\[
\{ k\} 2k (n-2k )\sigma^{ k+1 } [(2k)' + (n-2k)] (z) = 2k (n-2k
)\left(\sigma\tD_i^k (z) -\sum_{ j=0 }^{ k-1 } b_j^k\sigma^{ k+1+j
}\tD_i^j (z)\right).
\]
Inserting these into \eqref{eqnknevenrecur} we obtain
\begin{equation}\label{eqnk+1neven}
\begin{split}
\tD_i^{ k+1 } (z) = \left( -{ 2k\choose 2}{ n-2k+2\choose 2 } b_0^{
k-1} \sigma^{ k-1 } - 2k (n-2k) b_0^k\sigma^{ k+1 }\right) z \\ 
+\left(
a_0^k\sigma^k - {2k\choose 2}{ n-2k+2\choose 2} b_1^{ k-1 } \sigma^k -
2k (n-2k ) b_1^k\sigma^{ k+2 }\right)\tD_i (z) +\\
\vdots \\
+\left( b_{ k-3 }^k\sigma^{ 2k-3 } - {2k\choose 2}{ n-2k+2\choose 2}
b_{ k-2 }^{ k-1 }\sigma^{ 2k-3 } - 2k (n-2k ) b_{ k-2 }^k \sigma^{
2k-1 }\right)\tD_i^{ k-2 } (z) \\
+\left( b_{ k-2 }^k\sigma^{ 2k-2 } + {2k\choose 2}{ n-2k+2\choose 2} - 
2k (n-2k ) b_{ k-1 }^k\sigma^{ 2k }\right)\tD_i^{ k-1 } (z) \\
+ \left( b_{ k-1 }^k\sigma^{ 2k-1 } + 2k (n-2k )\sigma\right)\tD_i^k
(z) \\
+ \{ k+1\} \sigma^{ k+1 } [(2k+2)' + (n-2k-2) ] (z).
\end{split}
\end{equation}

For $k=2$ we only need to replace $6[ 4+ (n-4)'] (z)$ which is
\[
6 [ 4+ (n-4)'] (z) = \sigma\tD_i^2 - {n\choose 2}\sigma -2 (n-2 )\tD_i
\]
and we obtain the equation
\begin{equation*}
\begin{split}
\tD_i^3 = -4 (n-4 ){n\choose 2 }\sigma +\left({n\choose 2} + 6
{n-2\choose 2 } - 4(n-4 ) 2 (n-2 )\right)\tD_i \\
+ \left( 2 (n-2 ) +
4(n-4)\right)\sigma\tD_i^2 + 6\cdot 15 [6 + (n-6)'].
\end{split}
\end{equation*}
This proves the existence of (\ref{eqn3.3})
for all $k\leq\frac{n-2}{4}$.
\end{proof}

\subsection{The recursion relations between the coefficients for $n$
even}\label{ss2.2}
Using equations (\ref{eqnD2even}) and
(\ref{eqn3.3}) we obtain the following
initial values
\[
b_0^2 = {n\choose 2} \quad \mbox{and} \quad b_1^2 = 2 (n-2). 
\]
Using equations (\ref{eqn3.3}) and (\ref{eqnk+1neven}) we 
obtain, for  $2\leq k\leq\frac{n-2}{4}$, the recursion relations for
the integers $b^\ell_j$.
\begin{corollary}\label{corcoefeven}
Setting $b^k_{-1} =b^{k-1}_k=0$ and $b^k_k = -1$ for $1 \leq k \leq
\frac{n-2}{2}$, we have for all $0 \leq j \leq k$
\[
b_{ j }^{ k+1 } = b_{ j-1 }^k - {2k\choose 2}{
n-2k+2\choose 2} b_{ j }^{ k-1 } - 2k (n-2k ) b_{ j}^k 
\]
\end{corollary}

\subsection{The final equations for $n$ even}
Suppose first that $n=4k-2$, $k\geq 2$. Then we have
\[
[(2k)' + (n-2k)] = [(2k)' + (2k-2)] =\sigma [(2k-2)' + 2k ] = \sigma.
[(2k-2)' + (n-2k+2)]
\]
So, combining \eqref{eqn3.3} for $k$ and $k-1$, we obtain
\begin{proposition} \label{prop3.4}
Suppose $n=4k-2, k \geq 2$. Then $\tD^i$ satisfies the following equation
\begin{equation}\label{eqnn4k-2}
X^k - \sum_{j=0}^{k-1}\left( b^k_j \sigma^{k+j} - {2k \choose 2}b^{k-1}_j
\sigma^{k+j-1}\right) X^j = 0,
\end{equation}
where the $b^\ell_j$ are the integers of subsection \ref{ss2.2}.
\end{proposition}

Now suppose $n=4k$, $k\geq 2$.
Here we apply $\tD_i$ to \eqref{eqn3.3} for $k=\frac{n}{4}$ to
obtain
\[
\begin{array}{ll}
\tD_i^{ k+1 } &= \sum_{j=0}^{k-1} b_j^k\sigma^{k+j}\tD_i^{j+1} + \{
k\}\sigma^k\tD_i [(2k)' + 2k] \\
&= \sum_{j=0}^{k-1} b_j^k\sigma^{k+j}\tD_i^{j+1} +\\
&+\{ k\}\sigma^k\left( {2k+2\choose 2}\left( [ (2k+2)' +
(2k-2) ] + [ (2k-2)' + (2k+2) ]\right) + 4k^2 [ (2k)' + 2k ]\right).
\end{array}
\]
Now we use equation \eqref{eqn3.3} for $k-1$ and its image by
$\sigma$ to replace $[ (2k+2)' + (2k-2) ] + [ (2k-2)' + (2k+2) ]$ and
equation \eqref{eqn3.3} for $k$ to replace $[ (2k)' + 2k ]$ and
obtain
\begin{proposition} \label{prop3.5}
Suppose $n= 4k, k \geq 2$. Then $\tD_i$ satisfies the following equation
\begin{equation}\label{eqnn4k}
\begin{split}
X^{k+1} +
 \sum_{j=0}^k \left( { 2k\choose 2}{2k+2\choose 2 } b_j^{ k-1} -
b^k_{j-1} \right)\sigma^{k+j-1}X^j\\
  + \sum_{j=0}^k \left( { 2k\choose 2}{2k+2\choose 2 } b_j^{k-1 } +
4k^2 b_j^k\right)\sigma^{ k+j }X^j = 0. \\
\end{split}
\end{equation}
where the $b^\ell_j$ are the integers of subsection \ref{ss2.2}.
\end{proposition}

\subsection{The equations in $B^{\sigma}_i$ and $P^{\sigma}_i$.}
According to subsection \ref{ss3.2} the correspondences $\tD_i$ induce
endomorphisms on the abelian subvarieties
$B^{\sigma}_i = Im(\sigma+1) \subset J\tC_i$ and $P^{\sigma}_i =
Im(\sigma - 1) \subset J\tC_i$
which we denote by the same letter.  

On $B^{\sigma}_i$ we have $\sigma =1$. Inserting this into Propositions 
\ref{prop3.4} and \ref{prop3.5} we finally ontain the following result.

\begin{theorem} \label{theorem3.6}
On the abelian variety $B^{\sigma}_i$ the endomorphism $\tD_i$
satisfies the following equation,\\
{\em (1)} for $n=4k-2, k \geq 2$,
\begin{equation} \label{eqn3.8}
X^k + \sum_{j=0}^{k-1}\left(b^k_j - {2k \choose 2}b^{k-1}_j \right) X^j = 0,
\end{equation}
{\em (2)} for $n = 4k, k \geq 2$,
\begin{equation} \label{eqn3.9}
X^{k+1} + \sum_{j=0}^k \left(2 {2k \choose 2}{2k+2 \choose 2}
b^{k-1}_j + 4k^2b^k_j - b^k_{j-1} \right) X^j = 0.
\end{equation}
\end{theorem}

On $P^{\sigma}_i$ we have $\sigma = -1$. Here we obtain 
\begin{theorem} \label{theorem3.7}
On the abelian variety $P^{\sigma}_i$ the endomorphism $\tD_i$
satisfies the following equation,\\
{\em (1)} for $n=4k-2, k \geq 2$,
\begin{equation} \label{eqn3.10}
X^k + \sum_{j=0}^{k-1} (-1)^{k+j} \left(b^k_j + {2k \choose
2}b^{k-1}_j \right) X^j = 0,
\end{equation}
{\em (2)} for $n = 4k, k \geq 2$,
\begin{equation} \label{eqn3.11}
X^{k} - \sum_{j=0}^{k-1} b^k_j \tD_i^j = 0.
\end{equation}
\end{theorem}
Note that after proving Theorem \ref{mainthm2} we can conclude that equation \eqref{eqn3.11} means that the eigen-abelian
variety of one of the roots of equation
\eqref{eqnn4k} on $P_i^{\sigma}$ has dimension $0$.
\begin{proof}
(1) is a direct consequence of Proposition \ref{prop3.4}. For $n = 4k,
n \geq 2$ we obtain an equation of degree $k$ by noting that $[(2k)' +
2k] = \sigma[(2k)'+2k]$. Subtracting equation \eqref{eqn3.3} from its
own image by $\sigma$ and dividing by $-2$ we obtain the equation
\eqref{eqn3.11} on $P^{\sigma}_i$ after replacing $\sigma$ by $-1$.
\end{proof}

\section{Combinatorial preliminaries}

\noindent
In order to find the zeros of \eqref{eqn2.8}, \eqref{eqn2.9}, \eqref{eqn3.8}, \ldots,  \eqref{eqn3.11} we need some combinatorial   
properties relating our set up to the \textsc{Hamming} scheme 
from algebraic graph theory  (see \cite{macwilliams}, \cite{godsil} for background
information). In particular, we shall use the fact that the eigenvalues of the distance$-k$ transform are
given by values of the \textsc{Krawtchouk} polynomials.
For convenience, we will keep the presentation self-contained.
Note that for the proofs of Theorems 1 and 2 only the cases 
$k=n-1$ and $k=n-2$ below are relevant.

\subsection{The distance$-k$ transform and its eigenvalues}
Consider the group 
\[\mathbb{B}^n = \mathbb{Z}_2^n = (\{0,1\}^n,\oplus)
\] 
of bitvectors of length $n$ with componentwise addition mod 2. For 
$x =(x_1,\ldots,x_n) \in \mathbb{B}^n$ and $y = (y_1,\ldots,y_n) \in \mathbb{B}^n$ let 
\[
\|x\|=\sum_{1 \leq i \leq n} x_i \qquad \mbox{and} \quad d(x,y) = \|x-y\|
\]
denote their Hamming weight and distance. Let $\mathbb{B}_k^n$ denote the set of
bitvectors of length $n$ and weight $k$ where $0 \leq k \leq n$.

For any field $F$ (below we assume that the characteristic 
of $F$ is  $\neq 2$), 
let 
\[R_n=F[\mathbb{B}^n]
\]
denote the
vector space over $F$ with $\mathbb{B}^n$ as a basis.
We consider the following endomorphisms of $R_n$:
\begin{itemize}
\item
The {\it Hadamard transform} is the endomorphism of $R_n$ defined 
on basis elements $x \in \mathbb{B}^n$ by
\[
x \mapsto \widehat{x} = \sum_{y \in \mathbb{B}^n} (-1)^{x \cdot y} \, y
 = \sum_{y \in \mathtt{B}^n} \chi_x(y)\,y,
\]
where $x \cdot y$ is the scalar product, i.e.,
$\chi_x :y \mapsto (-1)^{x\cdot y}$ denotes the character of
$\mathbb{B}^n$ belonging to $x$.
\item For $0 \leq k \leq n$ the {\it distance$-k$ transform} $\Gamma_{n,k}$
is the endomorphism of $R_n$ defined on basis elements $x \in
\mathbb{B}^n$ by
\[
x \mapsto \Gamma_{n,k} (x) = \sum_{y \in \mathbb{B}_k^n} x \oplus y,
\]
In other words $\Gamma_{n,k}$ associates to $x \in \mathbb{B}^n$ the
sum of all basis elements
at Hamming distance $k$ from $x$ (changing $k$ coordinates from $0$ to $1$ or
vice versa).
\end{itemize}

\begin{proposition} \label{prop4.1}
For $0 \leq k \leq n$, and
for $x \in \mathbb{B}_\ell^n~~(0 \leq \ell \leq n)$, the Hadamard transform
$\widehat{x}$ is an eigenvector of $\Gamma_{n,k}$ with
eigenvalue
\[
\lambda_{n,k,\ell} = \chi_x\left( \mathbb{B}_k^n \right)=
\sum_i (-1)^i \binom{\ell}{i} \binom{n-\ell}{k-i}.
\]
\end{proposition}

\begin{proof} First note that all the operators $\Gamma_{n,k} ~(0 \leq
k \leq n)$ commute,
hence they have a common system of eigenvectors. Write
\begin{align*}
\Gamma_{n,k} (\widehat{x}) &=
\sum_{y \in \mathbb{B}^n} (-1)^{x \cdot y}\Gamma_{n,k} (y) =
\sum_{y \in \mathbb{B}^n} (-1)^{x \cdot y} \sum_{z \in \mathbb{B}_k^n}
y \oplus z =\cr
&= 
\sum_{y \in \mathbb{B}^n}  \sum_{z \in \mathbb{B}_k^n} (-1)^{x \cdot
(y \oplus z)}y =
\sum_{y \in \mathbb{B}^n} (-1)^{x \cdot y} 
\left( 
\sum_{z \in \mathbb{B}_k^n} (-1)^{x \cdot  z} \right) y =\cr
&= \chi_x(\mathbb{B}_k^n) \cdot \widehat{x},
\end{align*}
where we have used $(-1)^{x \cdot (y \oplus z)} = (-1)^{x\cdot
y}(-1)^{x \cdot z}$.

It is clear that the eigenvalue 
\[
\lambda_{n,k}(x) = \sum_{z \in \mathbb{B}_k^n} (-1)^{x \cdot  z}
= \chi_x\left( \mathbb{B}_k^n \right)
\]
corresponding to $x$ depends only
on the weight $\|x\|=\ell$ of $x$, so that one can
write $\lambda_{n,k,\ell}$ for it. Now for $x \in \mathbb{B}_\ell^n$:
\[
\sum_{k=0}^n \lambda_{n,k,\ell}\, t^k =
\sum_{z \in \mathbb{B}^n} (-1)^{x \cdot z}\, t^{\|z\|} 
= (1-t)^{\ell}(1+t)^{n-\ell},
\]
from which the above expression for $\lambda_{n,k,\ell}$ follows
by comparison of coefficients of $t^k$.
\end{proof}

\begin{remark} {\em For  $n \in \mathbb{N}$  
the \textsc{Krawtchouk} polynomials $P_k(x;n)$ $(0 \leq k \leq n)$ are
defined by
\[
\sum_{0 \leq k \leq n} P_k(i;n) \, z^k = (1-z)^i(1+z)^{n-i},
\]
or equivalently,
\[
P_k(x;n) = \sum_{j=0}^k (-1)^j \binom{x}{j}\binom{n-x}{k-j},
\]
so that $\lambda_{n,k,\ell}=P_k(\ell;n)$. We note the following well known and easily proved properties of these eigenvalues.
\begin{align*}
\lambda_{n,k,\ell}&= (-1)^k \cdot \lambda_{n,k,n-\ell} ,\cr
\lambda_{n, k,\ell} &= (-1)^\ell \cdot \lambda_{n,n-k,\ell} ,\cr
\binom{n}{\ell} \cdot \lambda_{n,k,\ell} &= \binom{n}{k} \cdot
\lambda_{n,\ell,k}.
\end{align*}}
\end{remark}

\subsection{$\mathcal{S}_n$-symmetry}

Since the Hadamard transform and the distance$-k$ transforms 
are compatible with the natural action of the symmetric group 
$\mathcal{S}_n$ on $\mathbb{B}_n$ and on $R_n$, one can take
quotients and consider the vector space 
$$
\widetilde{R}_n = R_n/\mathcal{S}_n.
$$
It is convenient to take a polynomial model for this space, i.e., 
let 
$$
\cH_n = \cH_n(X,Y)
$$ 
denote the vectorspace
of homogeneous polynomials in variables $X,Y$ of degree $n$.
Take the monomials 
$$
\xi_\ell = X^\ell Y^{n-\ell}, \qquad \qquad  0 \leq \ell \le n
$$ 
as a basis,
where $\xi_\ell$ is taken as the image of the elements of $\mathbb{B}_k^n$.
Then the quotient action of the distance$-k$ transform has the matrix
representation
$G_{n,k}= \left[ g_{\ell,i} \right]_{0 \leq \ell,i \leq n}$, where
\[
g_{\ell,i} = 
\begin{cases}
\binom{\ell}{j} \binom{n-\ell}{k-j} & \text{if}~i=k+\ell-2j, \cr
0 & \text{otherwise}.
\end{cases}
\]
The quotient action of the distance$-k$ transform may also be 
represented as a differential operator
\[
\Delta_{k} = \frac{1}{k!}\,\sum_{j=0}^k \binom{k}{j}
X^jY^{k-j}D_X^{k-j}D_Y^j
\]
on  $\cH_n$. Then the eigenvectors take the convenient form
\[
v_{n,\ell} = v_{n,\ell}(X,Y) = (X-Y)^\ell(X+Y)^{n - \ell} ~~(0 \leq
\ell \leq n).
\]
These polynomials $v_{n,\ell}$ form a basis of $\cH_n$ 
adapted to the operators $\Delta_k$,
with eigenvalues $\lambda_{n,k,\ell}~(0 \leq k,\ell \leq n)$.

\begin{remark} {\em In terms of the \textsc{Krawtchouk} polynomials
\[
v_{n,\ell} =\sum_{k=0}^\ell P_k(\ell;n) X^kY^{n-k} = \sum_{k=0}^n
\lambda_{n,k,\ell} X^k Y^{n-k}.
\]
The remarkable fact that the $\lambda_{n,k,\ell}$ appear both as
coefficients of the eigenpolynomials $v_{n,\ell}$ and as their
eigenvalues corresponding to $\Delta_k$:
\[
\Delta_k v_{n,\ell} = \lambda_{n,k,\ell} \cdot v_{n, \ell}~~(0 \leq
k,\ell \leq n),
\]
can be written equivalently as
\[
\sum_{\ell=0}^n \binom{n}{\ell} \lambda_{n,k,\ell} \lambda_{n,j,\ell}
= 2^n \binom{n}{k} \delta_{k,j},
\]
which is the orthogonality relation for the polynomials $P_k(x;n)~(0
\leq k \leq n)$.}
\end{remark}

\subsection{More symmetry}

Let $\cH_n^+$, respectively $\cH_n^-$, denote the subspace of
symmetric, respectively antisymmetric, polynomials in $\cH_n$,
i.e.
\[ 
\cH_n^{\pm} = \{ p \in \cH_n\;|\;   p(X,Y)= \pm p(Y,X) \}.
\]
Then obviously $\{v_{n,2\ell}\,;\,0 \leq 2\ell \leq n\}$ is a basis
of $\cH_n^+$, and $\{v_{n,2\ell+1}\,;\,0 \leq 2\ell+1 \leq n\}$ is a
basis of $\cH_n^-$. Since the operators $\Delta_k$ are symmetric with
respect to $X,Y$, the subspaces $\cH_n^+$ and $\cH_n^-$ are
$\Delta_k$-invariant.

 Let $\cH_n^e$, respectively $\cH_n^o$, denote the subspace of
 polynomials in $\cH_n$ where the variable $Y$ appears only with even,
 respectively odd, powers, i.e.
\[
\cH_n^e = \{ p \in \cH_n \;| \; p(X,Y)=p(X,-Y)\}, \quad \cH_n^o = \{ p
\in \cH_n \;| \; p(X,Y)=-p(X,-Y)\}.
\]
Let
\[
p^e(X,Y)= (p(X,Y)+p(X,-Y))/2 \quad \mbox{and} \quad p^o(X,Y)=
(p(X,Y)-p(X,-Y))/2
\] 
denote the even and odd part of $p(X,Y)$.
We have $v_{n,\ell}(X,-Y)=v_{n,n-\ell}(X,Y)$, hence
\begin{align*}
v_{n,\ell}^e  &= (v_{n,\ell}+v_{n,n-\ell})/2, \cr
v_{n,\ell}^o   &= (v_{n,\ell}-v_{n,n-\ell})/2,
\end{align*}
so that the 2-dimensional subspace spanned by
$\{v_{n,\ell},v_{n,n-\ell}\}$ has also $\{v_{n,\ell}^e,v_{n,\ell}^o\}$
as a basis. A degenerate situation occurs for $n$ even and $\ell=n/2$,
where $v_{n,n/2}$ itself is an even polynomial, hence
$v_{n,n/2}^e=v_{n,n/2}$ and $v_{n,n/2}^o=0$, and we only have a
one-dimensional subspace.

From $\lambda_{n,k,\ell}=(-1)^k \cdot \lambda_{n,k,n-\ell}$ it follows that
\[
\Delta_k v_{n,\ell}^e = 
\begin{cases}
\lambda_{n,k,\ell} \cdot v_{n,\ell}^e &\text{if $k$ is even},\cr
 \lambda_{n,k,\ell} \cdot v_{n,\ell}^o &\text{if $k$ is odd},
 \end{cases}
\]
and similarly 
\[
\Delta_k v_{n,\ell}^o = 
\begin{cases}
\lambda_{n,k,\ell} \cdot v_{n,\ell}^o &\text{if $k$ is even},\cr
 \lambda_{n,k,\ell} \cdot v_{n,\ell}^e &\text{if $k$ is odd}.
 \end{cases}
\]

Hence, if $k$ is odd, then the subspaces $\cH_n^e$ and $\cH_n^o$ are
not $\Delta_k$-invariant, they are rather  $\Delta_k^2$-invariant with
\[
\Delta_k^2 v_n^e = \lambda_{n,k,\ell}^2 \cdot v_n^e, \quad \mbox{and} \quad
\Delta_k^2 v_n^o = \lambda_{n,k,\ell}^2 \cdot v_n^o.
\]

Moreover, if $k$ is even and $n$ is odd, then $\{ v_{n,\ell}^e\,;\,0
\leq \ell <n/2\}$ is a basis of $\cH_n^e$, and $\{ v_{n,\ell}^o \,;\,0
\leq \ell <n/2\}$ is a basis of $\cH_n^o$;

Finally, if $k$ and $n$ are even, then a basis of $\cH_n^e$ is $\{
v_{n,\ell}^e\,;\,0 \leq \ell <n/2\} \cup\{ v_{n,n/2} \}$, and a basis
for $\cH_n^o$ is the same as for $n$ odd.  Note that in this case
$v_{n,2\ell}^e \in \cH_n^+$ and $v_{n,2\ell+1}^e \in \cH_n^-$, and,
similarly, $v_{n,2\ell}^o \in \cH_n^+$ and $v_{n,2\ell+1}^o \in
\cH_n^-$, because $n$ and $n - \ell$ have the same parity.
 
As a consequence we obtain
\begin{proposition} \label{prop4.4}
{\em (1)} If $n$ is odd and $k$ is even, then the actions of $\Delta_k$ on the four invariant subspaces
$\cH_n^+$, $\cH_n^-$, $\cH_n^e$, $\cH_n^o$ of dimension $(n+1)/2$
are isomorphic, as they all
afford the $\lambda_{n,k,\ell}$ with $0 \leq \ell <n/2$
as eigenvalues.\\
{\em (2)} If $n$ and $k$ are both even, then the invariant subspaces
$\cH_n^+$  and $\cH_n^-$ (of dimension $n/2$)
are not only different in dimension $(\dim \cH_n^+ =n/2+1,
\dim \cH_n^- =n/2)$, but the actions of
$\Delta_k$ on these subspaces have complementary subsets of eigenvalues:
$\{\lambda_{n,k,2 \ell}\,;\,0 \leq 2 \ell \leq n/2 \}$ for $\cH_n^+$ and 
$\{\lambda_{n,k,2 \ell+1}\,;\,0 \leq 2 \ell +1 \leq n/2 \}$ for $\cH_n^-$.
All of these are double eigenvalues, except $\lambda_{n,k,n/2}$,
which is simple. 
\end{proposition}
In the case of Proposition \ref{prop4.4} (2) one can use $\cH_n^e$,
$\cH_n^o$ to separate the eigenvalues as follows.
Consider the invariant
subspaces 
\[
\cH_n^{+e} = \cH_n^+ \cap \cH_n^{e},\quad 
\cH_n^{+o} = \cH_n^+ \cap \cH_n^{o},\quad
\cH_n^{-e} = \cH_n^- \cap \cH_n^{e},\quad
\cH_n^{-o} = \cH_n^- \cap \cH_n^{o}.
\]
Then we have the following.
\begin{itemize}
\item if $n \equiv 0 \bmod 4$: 
\[
\begin{array}{r|clll}
                   & \dim & \text{eigenvalues} & \text{eigenvectors}
&\text{range}\cr\hline
\cH_n^{+e} & n/4+1 &  \lambda_{n,k,2\ell}     &  v_{n,k,2\ell}^e  & 0
\leq \ell \leq n/4\cr
\cH_n^{+o} & n/4    &  \lambda_{n,k,2\ell}     &  v_{n,k,2\ell}^o  &
0 \leq \ell < n/4\cr
\cH_n^{-e} & n/4     &  \lambda_{n,k,2\ell+1}  &  v_{n,k,2\ell+1}^e &
0 \leq \ell < n/4 \cr
\cH_n^{-o} & n/4     &  \lambda_{n,k,2\ell+1}  &  v_{n,k,2\ell+1}^o  &
0 \leq \ell < n/4
\end{array}
\]
\item if $n \equiv 2 \bmod 4$:
\[
\begin{array}{r|clll}
                   & \dim & \text{eigenvalues} & \text{eigenvectors}
&\text{range}\cr\hline
\cH_n^{+e} & (n+2)/4   &  \lambda_{n,k,2\ell}     &  v_{n,k,2\ell}^e
& 0 \leq \ell  < n/4\cr
\cH_n^{+o} & (n+2)/4    &  \lambda_{n,k,2\ell}     &  v_{n,k,2\ell}^o
&   0 \leq \ell < n/4\cr
\cH_n^{-e} & (n+2)/4     &  \lambda_{n,k,2\ell+1}  &
v_{n,k,2\ell+1}^e &   0 \leq \ell < n/4 \cr
\cH_n^{-o} & (n-2)/4     &  \lambda_{n,k,2\ell+1}  &
v_{n,k,2\ell+1}^o  &   0 \leq \ell < (n-2)/4    
\end{array}
\]
\end{itemize}

\section{Proofs of the main theorems}

\subsection{Proof of Theorem \ref{mainthm1}}

Let the situation be as in Theorem \ref{mainthm1}, i.e., suppose $n =
2k+1 \geq 3$ and consider the coverings of smooth projective curves
$\tX\stackrel{\kappa}{\lra} X\stackrel{\rho_n}{\lra }Y$ with $\rho_n$
of degree n and $\kappa$ \'etale of degree 2, satisfying Hypothesis \ref{hyp}.  Let $f: C \ra Y$ denote
the associated covering of degree $2^{n-1}$. For a point $y \in Y$ let
$\rho_n^{-1} (y)= \{{\overline x}_1, \ldots, {\overline x}_n \}$ and
$\kappa^{-1}({\overline x}_i) = \{x_i,x'_i \}$.  Then

\begin{equation}\label{eqnfibre}
f^{-1}(y) = \{ x_1^{\epsilon_1} + \cdots + x_n^{\epsilon_n} \;|\;
\epsilon = \;'\, \mbox{or no}\; ', \mbox{even number of\; $'$s} \}.
\end{equation} 
Since the correspondence $D$ is independent of the point $y \in Y$, we
can identify $f^{-1}(y)$ with the set of bit vectors of length $n$
with an even number of components different from 0:

\begin{equation} \label{eqnbit}
f^{-1}(y) = \mathbb{B}^{n,e} := \{(e_1, \ldots, e_n) \;|\; e_i = 0 \;
\mbox{or} \; 1, \sum e_i \; \mbox{even} \}.
\end{equation}
This gives us two additional structures on $f^{-1} (y)$, namely the addition $\oplus$
and the Hamming distance on $\mathbb{B}^{n.e}$.  Denote by $R_n^e$ the
corresponding subspace of $R_n$:
\[
R_n^e := F[\mathbb{B}^{n,e}],
\]
and define $\cH_n^e$ as in Section 4. Using these identifications the
correspondence $D$ on $C$ induces the distance$-(n-1)$ transform
$\Gamma_{n,n-1}$ on the vector spaces $R_n^e$ as well as the
differential operator $\Delta_{n-1}$ on the vector space
$\cH_n^e$. Now, according to Theorem \ref{thm2.4}, the correspondence
$D$ satisfies an equation of degree $k$. Hence to complete the proof
of Theorem \ref{mainthm1} it suffices to show that $\Delta_{n-1}$
admits the $k$ eigenvalues $(-1)^{k+j}(2j+1), \; 0 \leq j \leq
k$. This follows immediately from Proposition \ref{prop4.4}, since the only nonzero terms of
\[
\lambda_{n,n-1,\ell} =
\textstyle 
\sum_{j} (-1)^j \binom{\ell}{j} \binom{n-\ell}{n-1-j}
\]
are those where $j=\ell-1$ and $j=\ell$
as nonzero summands,  so that
\begin{align*}
\lambda_{n,n-1,\ell} &= 
\textstyle
(-1)^\ell \left( \binom{\ell}{0}\binom{n-\ell}{1}-
\binom{\ell}{1}\binom{n-\ell}{0} \right)\cr
&=\textstyle (-1)^\ell (n-2 \ell)  \cr
&=(-1)^{k+\ell}(2(k-\ell) +1) \quad \mbox{for} \quad0 \leq \ell \leq n.
\end{align*}
Setting $i=k-l$, this finishes the proof.
\hfill $\Box$

\subsection{Proof of Theorem \ref{mainthm2}}

Let the situation be as in Theorem \ref{mainthm2}, i.e., suppose $n = 2k
\geq 4$ and consider the coverings of smooth projective curves
$\tX\stackrel{\kappa}{\lra} X\stackrel{\rho_n}{\lra }Y$ with $\rho_n$
of degree n and $\kappa$ \'etale of degree 2, satisfying Hypothesis \ref{hyp}.  For $i=1,2$ let $f_i:
\tC_i \ra Y$ denote the associated covering of degree $2^{n-1}$. For a
point $y \in Y$ the fibre $f_i^{-1}(y)$ is given as in
\eqref{eqnfibre} and will be identified with $\mathbb{B}^{n,e}$ as in
\eqref{eqnbit}.  Defining $R_n^e$ and $\cH_n^e$ as in Section 5.1
and using these identifications, the correpondence $\tD_i$ on $\tC_i$
induces the distance$-(n-2)$ transform $\Gamma_{n,n-2}$ on the vector
space $R_n^e$ as well as the differential operator $\Delta_{n-2}$ on
the vector space $\cH_n^e$. Since $\sigma = 1$ on the abelian
subvariety $B_i^{\sigma}$ and $\sigma = -1$ on $P_i^{\sigma}$,
this implies that under the assumption $\sigma = 1$ the correspondence
$\tD_i |_{ B^{\sigma }_i}$ induces the operator $\Delta_{n-2}$ on the subspace
$\cH_n^{+e}$, and, similarly, under the assumption $\sigma =-1$ the
correspondence $\tD_i|_{P_i^{\sigma}}$ induces the operator
$\Delta_{n-2}$ on the subspace $\cH_n^{-e}$.

Suppose first that $n = 4k \geq 4$. According to Theorems
\ref{theorem3.6} and \ref{theorem3.7}, $\tD_i|_{B_i^{\sigma}}$,
respectively $\tD_i|_{P_i^{\sigma}}$, satisfies an equation of degree
$k+1$, respectively $k$. Hence it suffices to show that $\Delta_{n-2}$
admits the $k+1$ distinct eigenvalues $8(k- \ell)^2 - 2k, \; 0 \leq l \leq k$
on the vector space $\cH_n^{+e}$ and the $k$ distinct eigenvalues $-8(k-\ell)^2
+10k -8\ell -2, 0 \leq \ell \leq k-1$ on the vector space
$\cH_n^{-e}$.

According to the table at the end of Section 4, the eigenvalues of
$\Delta_{n-2}$ are $\lambda_{n,n-2,2\ell}, 0 \leq \ell \leq k$, on
$\cH_n^{+e}$, and $\lambda_{n,n-2,2\ell+1}, 0 \leq \ell \leq k-1$, on
$\cH_n^{-e}$.\\ In the formula
\begin{equation} \label{eqnlambda}
\lambda_{n,n-2,\ell} =\textstyle 
\sum_{j} (-1)^j \binom{\ell}{j} \binom{n-\ell}{n-2-j}
\end{equation}
only the three terms for $j=\ell-2,\ell-1,\ell$
are nonzero. So, for $0 \leq \ell \leq n$, we obtain
\begin{align*}
\lambda_{n,n-2,\ell} &= \textstyle (-1)^\ell
\left(
\binom{\ell}{2}\binom{n-\ell}{0}
- \binom{\ell}{1}\binom{n-\ell}{1}
+\binom{\ell}{0}\binom{n-\ell}{2}\right) \cr
&=\textstyle
(-1)^{\ell} \left({ \binom{n-2\ell}{2}}-\ell \right).
\end{align*}
Hence 
\[
\lambda_{n,n-2,2\ell} = \textstyle \binom{4k-4\ell}{2} -2l = 8(k-
\ell)^2 - 2k
\]
\[
\lambda_{n,n-2,2\ell +1} = -(\textstyle \binom{4k-4\ell - 2}{2} -2l
-1) = -8(k-\ell)^2 +10k -8\ell -2
\]
which completes the proof for $n=4k$.\\ The proof for $n = 4k-2$ is
essentially the same. We have to compute the eigenvalues of
$\Delta_{n-2}$ on $\cH_n^{-e}$ and $\cH_n^{-e}$. According to the last
table in Section 4 and \eqref{eqnlambda} they are
\[
\lambda_{n,n-2,2\ell} = \textstyle \binom{4k-2-4\ell}{2} -2l = 8(k-
\ell)^2 - 10k +8l + 3
\]
\[
\lambda_{n,n-2,2\ell +1} = -(\textstyle \binom{4k-4\ell - 4}{2} -2l
-1) = -8(k-\ell)^2 +18k -16\ell -9.
\]
This completes the proof of Theorem \ref{mainthm2}. \hfill $\Box$

\subsection{The correspondences associated to the distance$-k$ transform for $k\leq n-3$}

For $k\leq n-3$ the associated correspondences are,
\begin{itemize}
\item when $n$ is odd,
\[
\tD_{i,k} :=\{ (x_1 +\ldots + x_n , x_1' +\ldots + x_k' + x_{ k+1 }
+\ldots + x_n )\}\subset \tC_i\times\tC_i
\]
with image $D_k$ in $C\times C$. Then the equations in
Proposition \ref{propDk} show that the eigenvalues of the associated
endomorphism of $JC$ can be computed from those of $D$ and the
eigen-abelian varieties are the same as those of $D$;
\item when $k$ and $n$ are even,
\[
\tD_{i,k} :=\{ (x_1 +\ldots + x_n , x_1' +\ldots + x_k' + x_{ k+1 }
+\ldots + x_n )\}\subset \tC_i\times\tC_i
\]
with reduced image $D_i$ in $C_i\times C_i$. Then the equations in
Proposition \ref{propDkeven} show that the eigenvalues of the
associated endomorphisms of $J\tC_i$ can be computed from those of
$\tD_i$ and the eigen-abelian varieties are the same as those of
$\tD_i$.

Note that odd values of $k$ will give us correspondences between
$\tC_i$ and $\tC_{3-i}$.
\end{itemize}
So, using Propositions \ref{propDk} and \ref{propDkeven} and the
calculations following them, we obtain yet more combinatorial identities.

\section{The dimensions of the eigen-abelian varieties}

Let $n \geq 3$ be an integer.  In order to have a unified statement, we consider, for odd $n$, the curve $C$
as $\tC_i$ and the correspondence $D$ as $\tD_i$ (see diagram \ref{diag1.3}). 
  We will write systems of linear equations whose solutions are the
dimensions of the eigen-abelian varieties of $\tD_i$. These equations
will be obtained by computing the analytic traces of the powers of $\tD_i$ in two
different ways. Finally, we use these equations to compute the dimensions for $n \leq 10$. 
First we see that the eigen-abelian varieties can be parametrized explicitly.

\subsection{\bf Geometric description of eigen-abelian varieties} 

Choose a point $y$ of $Y$ where $\rho_n$ is not branched and, let
$\ox_1,\ldots ,\ox_n$ be the points of $\rho_n^{ -1 }(y)$ and $x_1 ,
x_1' ,\ldots , x_n , x_n'$ their inverse images in $\tX$. Fix
$i=1$ or $2$ and assume that $x_1 +\ldots + x_n\in \tC_i$. As in
\eqref{eqnfibre} we identify the fiber $\mu^{-1} (y)$ when $n$ is odd,
resp. $(\mu_i\tau_i)^{-1}(y)$ when $n$ is even, with the set 
${\mathbb B}^{n,e}$ of bit
vectors of length $n$ with an even number of components different from
$0$. Conversely, let $t_{(e_1, \ldots, e_n)} = x_1^{e_1} + \cdots + x_n^{e_n}$ denote 
the point of $C$ corresponding to $(e_1, \ldots e_n) \in {\mathbb B}^{n.e}$ where $x_i^{ e_i } = x_i$ if $e_i = 0$ and $x_i^{ e_i } = x_i'$ if $e_i =1$.
 
As in Section 5, the correspondence $\tD_i$ on the curve $\tC_i$ induces the distance$-k$ transform 
$\Gamma_{n,k}$ on the vector space $R_n^e$ where $k=n-1$ if $n$ is odd and $k= n-2$ if $n$ is even. 

According to Proposition \ref{prop4.1}, for each $\ell \in \{0, \ldots, n\}$ and $x \in \mathbb{B}^n_{\ell}$
the Hadamard transform
\[
x \mapsto \widehat{x} = \sum_{y \in \mathbb{B}^n} (-1)^{x
\cdot y} \, y,
\]
is an eigenvector for $\Gamma_{n,k}$ with eigenvalue
\[
\lambda_{n,k,\ell} = \chi_{x}\left( \mathbb{B}_k^n \right)=
\sum_i (-1)^i \binom{\ell}{i} \binom{n-\ell}{k-i}.
\]

Under the identification of ${\mathbb B}^{n,e}$ with $\mu^{-1} (y)$,
resp. $(\mu_i\tau_i)^{-1}(y)$, the Hadamard transform $\widehat{x}$
corresponds to the divisor
\[
\sum_{y \in \mathbb{B}^n} (-1)^{x
\cdot y} \, t_{y}
\]
on $\tC_i$. Recall that $k=n-1$ if $n$ is odd and $k=n-2$ if $n$ is
even. We have proved

\begin{proposition}
The eigen-abelian subvariety of $J\tC_i$
for the eigenvalue $\lambda_{n,k,\ell}$ is generated by divisors of the
form $\sum_{y \in \mathbb{B}^n} (-1)^{x
\cdot y} \, t_{y}$ where $x \in {\mathbb B}^n_\ell$ is fixed, after
substracting a fixed divisor of the correct degree.
\end{proposition}

The map
\[
\begin{array}{rcl}
\tX & \inj & Div^{0}\tC_i \\
p & \mapsto & \left( p + \tX^{ (n-1 )}\right)\cap \tC_i - \left(
\sigma p + \tX^{ (n-1 )}\right)\cap \tC_i
\end{array}
\]
induces a map from the Prym variety $P(\tX\ra X)$ of $\tX\ra X$ to $J\tC$ which is easily seen to be an isogeny to its image. For $n$ odd, let $d_{\lambda }$ be the dimension of the eigen-abelian variety of
the eigenvalue $\lambda$ of $\tD_1$. For $n$ even, let $d_{\lambda }$, resp. $e_{\lambda }$, be the dimension of the
eigen-abelian variety of the eigenvalue $\lambda$ of $\tD_i$ in
$P_i^{\sigma }$, resp. $B_i^{\sigma}$.
For the special values $\ell=0$ and $\ell=1$ we obtain

\begin{corollary}\label{cor2.2}
1. For $n$ odd the eigen-abelian subvariety of $J\tC_i$ for the
eigenvalue $n$ is $\mu^* JY$. The eigen-abelian subvariety of $J\tC_i$
for the eigenvalue $-n+2$ is the image of $P(\tX\ra
X)$. In particular, $d_n = g_Y$ and $d_{-n+2} = g_X -1$.

2. If $n$ is even the eigen-abelian subvariety of $J\tC_i$ for the
eigenvalue ${n\choose 2}$ is $(\mu_i\tau_i)^* JY$. The eigen-abelian
subvariety of $J\tC_i$ for the eigenvalue $-\frac{ (n-1) (n-4) }{2}$
is the image of $P(\tX\ra X)$. In particular, $e_{n\choose 2} = g_Y$ and $d_{-\frac{(n-1)(n-4)}{2}} = g_X -1$.
\end{corollary}

\subsection{The case $n=2k+1$ odd}

Recall that $d_{\lambda }$ is the dimension of the eigen-abelian variety of
the eigenvalue $\lambda$ of $\tD_1$. Our first equation is
\[
\sum_{j=0 }^k d_{ (-1)^{j+k} (2j+1)} = g_{\tC_1}.
\]
Since we already know $d_n$ and $d_{-n +2}$ from Corollary
\ref{cor2.2}, we now need $k-2$ independent linear equations.

An endomorphism $D$ of an abelian variety $A$ naturally acts on the tangent space $T_0A$ of $A$ at the origin as well as on $H^1 (A ,\bQ)$. We denote $\tr_a (D)$ the analytic trace of $D$, i.e., the trace of $D$ as an endomorphism of $T_0A$, and $\tr_r(D)$ the rational trace of $D$, i.e., the trace of $D$ as an endomorphism of $H^1 (A ,\bQ)$.
Then, for every $\ell$ ($\leq k-2$), we have
\[
\tr_a (\tD_1^{\ell}) = \sum_{j=0 }^k ((-1)^{j+k} (2j+1))^{\ell} d_{
(-1)^{j+k} (2j+1)}.
\]
Now
\[
\begin{array}{c}
\tr_a (\tD_1^{\ell}) = \frac{1}{2} \tr_r (\tD_1^{\ell} ) = deg (\tD_1^{\ell}
) -\frac{1}{2}\Delta_{\tC_1 }\cdot \tD_1^{\ell}
\end{array}
\]
by \cite{bl04} page 334 Proposition 11.5.2. Since $deg
(\tD_1^{\ell} )= (deg (\tD_1 ))^{\ell}= n^{ \ell }$, it remains to
compute the intersection number $\Delta_{\tC_1 }\cdot \tD_1^{\ell}$ in
order to obtain a complete system of equations. For this we use
induction on $\ell$ and the equations of Proposition \ref{propDk}, namely
\begin{equation}\label{eqnelleven}
\ell \: even \hskip40pt \tD_1^{\ell } (z) = a_0^{\ell } z + a_2^{\ell }
\tD_1^2 (z) +\ldots + a_{ \ell-2 }^{\ell }\tD_1^{ \ell-2 }
(z) + {\ell}! [{\ell}' + (n-{\ell} )] (z)
\end{equation}

\begin{equation}\label{eqnellodd}
\ell \: odd \hskip40pt \tD_1^{\ell } (z) = a_1^{\ell} \tD_1 (z) +
a_3^{\ell} \tD_1^3 (z) +\ldots + a_{ {\ell}-2 }^{\ell}\tD_1^{ {\ell}-2
} (z) + {\ell}! [{\ell} + (n-{\ell} )'] (z)
\end{equation}
where the coefficients $a_{ {\ell}- 2m }^{\ell}$ are as in Proposition \ref{prop2.3}.

To compute the trace of $[{\ell}' + (n-{\ell} )]$, we need to count
the points $z = x_1 +\ldots + x_n $ such that, after possibly
renumbering the $x_i$, we have
\[
z = x_1' + \ldots + x_{\ell }' + x_{\ell +1} + \ldots + x_n.
\]
This happens only when $\ell =2$ and $x_2 = x_1'$, so that $\kappa
(x_1 )$ is a ramification point of $\rho_n$. Each such ramification
point gives $2^{ n-3 }$ points of $[{2}' + (n-2 )]\cdot\Delta_{\tC_1
}$. So we obtain
\[
[{\ell}' + (n-{\ell} )]\cdot\Delta_{\tC_1 } = 0
\]
for $\ell\neq 2$ and
\[
[{2}' + (n-2 )]\cdot\Delta_{\tC_1 } = 2^{ n-3 }\cdot deg (R_{ X/Y }) =
2^{ n-3 } (2 g_X -2 - 4( 2g_Y -2 )).
\]

\subsection{The case $n$ even}

We will treat the case $n=4k$, the case $n=4k-2$ is similar. Recall
that $d_{\lambda }$, resp. $e_{\lambda }$, is the dimension of the
eigen-abelian variety of the eigenvalue $\lambda$ of $\tD_i$ in
$P_i^{\sigma }$, resp. $B_i^{\sigma}$. Here we first have the two
equations
\[
\sum_{j=0 }^{k-1} d_{ -8 (k-j )^2 +10 k -8j-2} = g_{C_i} -1,
\]
\[
\sum_{j=0 }^{k} e_{ 8 (k-j )^2 -2k} = g_{C_i}.
\]
Since we already know $d_{-\frac{(n-1)(n-4)}{2}}$ and $e_{n\choose 2}$
from Corollary \ref{cor2.2}, we now need $2k-1$ independent linear
equations.

Here we compute the analytic traces of $D_i$ and $\tD_i$ in two
different ways. For every $\ell \leq k - 1$ we have 
\[
\tr_a (D_i^{\ell}) = \sum_{j=0}^k (8(k-j)^2 -2k)^{\ell}\; e_{8(k-j)^2
-2k}
\]
and
\[
\tr_a(\tD_i^{\ell}) = \sum_{j=0}^{k-1} (-8(k-1)^2 + 10k -8j
-2)^{\ell} \, d_{-8(k-1)^2 + 10k -8j
-2}+ \sum_{j=0}^k (8(k-j)^2 -2k)^{\ell}\, e_{8(k-j)^2
-2k}.
\]
On the other hand
\[
\begin{array}{c}
\tr_a (D_i^{\ell})  = \frac{1}{2} \tr_r (D_i^{\ell} ) = deg (D_i^{\ell}
) -\frac{1}{2}\Delta_{C_i }\cdot D_i^{\ell}
\end{array}
\]
and
\[
\begin{array}{c}
\tr_a (\tD_i^{\ell}) = \frac{1}{2} \tr_r (\tD_i^{\ell} ) = deg (\tD_i^{\ell}
) -\frac{1}{2}\Delta_{\tC_i }\cdot \tD_i^{\ell}
\end{array}
\]
by \cite[page 334 Proposition 11.5.2]{bl04}.  Since $deg
(D_i^{\ell} )= (deg (D_i ))^{\ell}= {n \choose 2}^{ \ell }$ and $deg
(\tD_i^{\ell} )= (deg (\tD_i ))^{\ell}={n \choose 2}^{ \ell }$, it remains to
compute the intersection numbers $\Delta_{C_i }\cdot D_i^{\ell}$ and $\Delta_{\tC_i }\cdot \tD_i^{\ell}$ in
order to obtain a complete system of equations. For this we use
induction on $\ell$ and the equations of Proposition \ref{propDkeven} which are in this situation
\begin{equation}\label{eqnkneven}
\tD_i^{\ell} (z) = \sum_{j=0}^{{\ell}-1}b_j^{\ell}\sigma^{{\ell}+j}\tD_i^j (z) 
+ \{ {\ell}\}\sigma^{\ell} [2{\ell}' + (n-2{\ell} )] (z).
\end{equation}
On $JC_i$ they become 
\begin{equation}
D_i^{\ell} (z) = \sum_{j=0}^{{\ell}-1}b_j^{\ell}\tD_i^j (z) 
+ \{ {\ell}\} [2{\ell}' + (n-2{\ell} )] (z).
\end{equation}
Here the coefficients $b_j^{\ell}$ are given by Corollary \ref{corcoefeven}.
The trace of $[2{\ell}' + (n-2{\ell} )]$ is the same as in the case
$n$ odd.

\subsection{The case $n=3$}

For $n=3$ the correspondence $\tD_1$ has two eigenvalues: $-1$
and $3$. By Corollary \ref{cor2.2} the eigen-abelian varieties are the images of the
Prym variety of $\tX \ra X$ and of $JY$ respectively.

\subsection{The case $n=4$}

For $n=4$ the correspondence $\tD_i$ has three eigenvalues: $0$ on
$P_i^{\sigma}$, $-2$ and $6$ on $B_i^{\sigma}$. 
By Corollary \ref{cor2.2}, the eigen-abelian variety for $0$ is the
image of the Prym $P(\tX\ra X)$ and the eigen-abelian variety for $6$ is
the image of $JY$. To compute $e_{-2}$ we use the equation 
\[
e_{-2} + e_6 = g_{C_i} = g_X
\]
from which it follows that $e_{-2} = g_X - g_Y$. 

In particular, in this case the three Prym varieties $P_1^{\sigma},
P_2^{\sigma}$ and $P(\tX \ra X)$ are isogenous.

\subsection{The case $n=5$}

When $n=5$ the correspondence $\tD_1$ has the three eigenvalues $-3,1,5$.
From Corollary \ref{cor2.2} we deduce that $d_{-3} = g_X -1$ and $d_5
= g_Y$. Now we have the equation
\[
d_{-3} + d_1 + d_5 = g_{\tC_1 } = 4 (g_X - g_Y ) +1
\]
which gives $d_1 = 3 (g_X -1) - 5 (g_Y-1)$.

\subsection{The case $n=6$}

For $n=6$ the correspondence $\tD_i$ has the eigenvalues $-5, 3$
on $P_i^{\sigma}$ and the eigenvalues $-1, 15$ on $B_i^{\sigma}$. We
already know that $d_{-5} = g_X - 1$ and $e_{15} = g_Y$. In addition
we have the two equations 
\[
\begin{array}{c}
d_{-5} + d_{3} = g_{ C_i } -1 = 4 ( g_X -1 -2 (g_Y -1 )) \\
e_{-1} + e_{15} = g_{ C_i } = 4 ( g_X -1 -2 (g_Y -1 )) +1
\end{array}
\]
which give
\[
d_3 = 3 (g_X -1 ) - 8 (g_Y -1 ),\quad
e_{-1} = 4 (g_X -1 ) - 9 (g_Y -1 ).
\]

\subsection{The case $n=7$}

For $n = 7$ the eigenvalues of $\tD_1$ are $-5,-1,3,7$. We know
$d_{-5} = g_X - 1$ and $d_7 = g_Y$. To compute $d_{-1}$ and $d_{-3}$ we
use the equations
\[
\begin{array}{c}
d_{-5} + d_{-1} + d_3 + d_7= g_{\tC_1 } = 16(g_X -1 -3(g_Y-1)) +1 \\
-5 d_{ -5 }  -  d_{ -1}  +  3 d_3  +  7 d_7  =  7\\
\end{array}
\]
to obtain
\[
d_{-1} = 10 (g_X -1) - 35 (g_Y-1), \quad d_3 =
5 (g_X -1) - 14 (g_Y-1).
\]

\subsection{The case $n=8$}

For $n=8$ the endomorphism $\tD_i$ has the eigenvalues $-14,2$ on
$P_i^{\sigma}$ and $-4,4,28$ on $B_i^{\sigma}$. We know $d_{-14} = g_X
- 1$ and $e_{28} = g_Y$. Here the additional equations are
\[
\begin{array}{rcl}
d_{-16} + d_{-14} + d_2 + e_{ -4 } + e_4 + e_{ 28 } & = & g_{ \tC_i } =
2^5 ( g_X -1 -4 (g_Y -1 )) +1 \\
e_{-4} + e_{4} + e_{ 28 } & = & g_{ C_i } = 2^4 ( g_X -1 -4 (g_Y -1
))+ 1 \\
- 16 d_{-16} - 14 d_{-14} + 2 d_2 - 4 e_{-4} + 4 e_4
+ 28 e_{ 28 } & = & 28, \\
\end{array}
\]
and we obtain
\[ 
d_2 = 15 (g_X -1) - 64 (g_Y-1),\quad 
e_{ -4 } = 10 (g_X -1) - 45 (g_Y-1),\quad   e_4 = 6 (g_X -1) - 20 (g_Y-1).
\]

\subsection{The case $n=9$}

For $n=9$ the eigenvalues of $\tD_1$ are $-7,-3,1,5,9$. We have
$d_{-7} = g_X -1$ and $d_9 = g_Y$. We have the equations
\[
\begin{array}{ccccccccccl}
d_{-7} & + & d_{-3} & + & d_1 & + & d_5 & + & d_9 & = & 2^6((g_X -1)
-5(g_Y-1)) +1 \\
-7 d_{ -7 } & - & 3 d_{ -3} & + & d_1 & + & 5 d_5 & + & 9 d_9 & = & 9\\
49 d_{ -7} & + & 9 d_{-3} & + & d_1 & + & 25 d_5 & + & 81 d_9 & = &
7\cdot 2^6 (g_X -1 ) - 3\cdot 9\cdot 2^6 (g_Y -1) + 81\\
\end{array}
\]
which give
\[
d_{-3} = 21 (g_X -1) - 105 (g_Y-1), \quad 
d_1 = 35 (g_X -1) - 189 (g_Y-1), \quad  d_5 = 7 (g_X -1) - 27 (g_Y-1).
\]

\subsection{The case $n=10$}

Here the eigenvalues of $\tD_i$ are $-27, -3, 5$ on $P_i^{\sigma}$, 
$-3, 13, 45$ on $B_i^{\sigma}$ and $d_{-27} = g_X -1$, $e_{45} = g_Y$.
The equations are
\[
\begin{array}{rrrrrrcl}
d_{-27} & + d_{-3} & + d_5 & + e_{-3} & + e_{ 13 } & + e_{ 45 }
& = & 2^7((g_X -1) -6(g_Y-1)) +1 \\
& & & e_{-3} & + e_{ 13 } & + e_{ 45 }
& = & 2^6((g_X -1) -6(g_Y-1)) +1 \\
- 27 d_{-27} & - 3 d_{-3} & + 5 d_5 & - 3 e_{-3} & + 13 e_{ 13 
}& + 45 e_{ 45 } & = & 45 \\
& & & -3 e_{-3} & + 13 e_{ 13 } & + 45 e_{ 45 }
& = & 45 - 2^6((g_X -1) -10(g_Y-1)) \\
\end{array}
\]
and we have
\[
\begin{array}{ll}
d_{ -3} = 28(g_X - 1) - 160 (g_Y -1), & d_5 = 35
(g_X -1) - 224 (g_Y-1), \\
e_{ -3 } = 56 (g_X -1) - 350 (g_Y-1), & e_{ 13 } = 8 (g_X -1) - 35 (g_Y-1).
\end{array}
\]

\section{A remark for the distance$-(n-1)$ transform for $n$ odd}
 
 In this particular case the action of $\Gamma_{n,n-1}$ on $\cH_n^+$
 (or on any of $\cH_n^-$ or $\cH_n^e$ or $\cH_n^o$)
 has the characteristic polynomial
\[
 C_n^+(X) = \prod_{\ell=1}^m \left(X-(-1)^{\ell+m}(2 \ell -1)\right),
\]
 where $m=(n+1)/2$. In terms of this action this means that,
\[
 C_n^+(X) = \det 
 \begin{bmatrix}
 X & -1 & 0 & 0 & \ldots & 0 &0 & 0\cr
 -n & X & -2 & 0 &  \ldots & 0 & 0& 0\cr
 0 & -n+1 & X & - 3 & \ldots& 0 & 0 & 0\cr
 \vdots & \vdots &\vdots & \vdots & \ddots & \vdots & \vdots & \vdots \cr
 0 & 0 & 0 & 0 & \ldots & X &  -m+2 & 0 \cr
 0 & 0 & 0 & 0 & \ldots & -m-2 & X & -m+1 \cr
 0 & 0 & 0 & 0 & \ldots & 0 &-m-1 &X-m
 \end{bmatrix}.
\]
 Note that the eigenvectors $\{ \lambda_{n,n-1,2 \ell}\,;\, 0 \leq
 \ell < m \}$ span the space under consideration. The fact that this
 determinant factors as mentioned before can be proved directly using
 elementary row and column operations with induction. On the other
 hand, there is a standard evaluation of the determinant of a
 tridiagonal matrix: let $a_1,a_2,\ldots$ and $b_1,b_2,\ldots$ be
 variables and define, for $m \geq 0$, the $(m+1)\times(m+1)$ matrix
 $M^{(m)}=M^{(m)}(a_1,a_2,\ldots;b_1,b_2,\ldots)$ by
\[
 M^{(m)} = 
 \begin{bmatrix}
 X & a_1 & 0 & 0 & \ldots & 0 & 0 & 0 \cr
 b_1 & X & a_2 & 0 & \ldots & 0 & 0 & 0 \cr
 0 & b_2 & X  & a_3 & \ldots & 0 & 0& 0 \cr
 \vdots & \vdots & \vdots & \ddots & \vdots & \vdots &\vdots \cr
 0 & 0 & 0 & 0 & \ldots & X & a_{m-1} & 0 \cr
 0 & 0 & 0 & 0 & \ldots & b_{m-1} & X & a_{m}\cr
 0 & 0 & 0 & 0 & \ldots & 0 & b_{m} & X
 \end{bmatrix}.
\]
Then
\[
\det M^{(m+1)} = X \cdot \det M^{(m)} - a_mb_m \det M^{(m-1)},
\]
which by induction gives
\[
 \det
  M^{(m)} = \sum_{0 \leq 2j \leq m} (-1)^j c_j^{(m)} X^{m-2j},
\]
 where
\[
 c_j^{(m)} = \sum_{ 1 \leq i_1 \ll i_2 \ll \cdots \ll i_j \leq m}
 a_{i_1}b_{i_1} a_{i_2}b_{i_2} \cdots a_{i_j}b_{i_j},
\]
 the notation $x \ll y$ meaning that $ x+1<y $. We note in passing that
 this sum has a combinatorial interpretation in terms of ``matchings'', 
 and the determinant evaluation as a second order recurrence points
 to a relation to orthogonal polynomials (see e.g. \cite{godsil}). 
 
Now
\[
C_n^+(X) = (X-m) \cdot \det M^{(m-1)} - (m^2-1) \cdot \det M^{(m-2)},
\]
where $(a_1,a_2,a_3,\ldots) = (-1,-2,-3,\ldots)$ and
$(b_1,b_2,b_3,\ldots) = (-n,-n+1,-n+2,\ldots)$.
 
 The fact that $C_n^+(X)$ factors into very simple linear factors leads
 to a surprisingly simple recurrence for the coefficients $c_j^{(m)}$
 in this particular situation, 
 which is not at all obvious from its definition, nor from its
 combinatorial interpretation.


\providecommand{\bysame}{\leavevmode\hbox to3em{\hrulefill}\thinspace}
\providecommand{\MR}{\relax\ifhmode\unskip\space\fi MR }
\providecommand{\MRhref}[2]{%
  \href{http://www.ams.org/mathscinet-getitem?mr=#1}{#2}
}
\providecommand{\href}[2]{#2}

\end{document}